\documentclass[12pt,leqno]{amsart}
\usepackage{amsfonts,amsmath,amsthm,amssymb,bbm,eufrak}
\usepackage[usenames,dvipsnames]{color}
\usepackage{enumerate,epsfig,graphics,indentfirst}
\usepackage{hyperref,latexsym,multicol}
\usepackage{pxfonts,rotating,yfonts}
\usepackage[labelsep=space]{caption}
\usepackage[all]{xy}
\usepackage{tikz}
\usetikzlibrary{decorations.markings}
\usepackage[labelformat=empty]{caption}
\usepackage{mathtools}

\newtheorem{theorem}{Theorem}[section]
\newtheorem{lemma}[theorem]{Lemma}

\theoremstyle{definition}

\usepackage[ansinew]{inputenc}
\usepackage{ifthen}
\usepackage{animate}

\newcommand{\be}{\begin{equation}}
\newcommand{\ee}{\end{equation}}

\begin{document}

\title[Diagonalizable Thue Equations]{Diagonalizable Thue Equations -- revisited}

\author[Saradha]{N. Saradha}
\address[1]{INSA Senior Scientist\\
DAE Centre for Excellence in Basic Sciences\\
University Of Mumbai,
Mumbai-400098, INDIA\\
Mailing address:\\
B-706, Everard Towers\\
Eastern Express Highway\\
Sion, Mumbai 400\,022 \textsc{India}}
\email{saradha54@gmail.com}

\author[Sharma]{Divyum Sharma}
\address[2]{Department of Mathematics\\
Birla~Inst{i}tute~of~Technology~and~Science, Pilani 333\,031 \textsc{India}}
\email{\texttt{divyum.sharma\symbol{64}pilani.bits-pilani.ac.in}}

\date{\today}
\subjclass[2010]{Primary 11D61}
\keywords{Diagonalizable forms, Thue equations}

\begin{abstract}
Let $r,h\in\mathbb{N}$ with $r\geq 7$ and let
$F(x,y)\in \mathbb{Z}[x ,y]$ be a binary form such that
    \[
    F(x , y) =(\alpha x + \beta y)^r -(\gamma x + \delta y)^r,
    \]
where $\alpha$, $\beta$, $\gamma$ and $\delta$  are algebraic constants with
$\alpha\delta-\beta\gamma \neq 0$.
We establish upper bounds for the number of primitive 
solutions to the Thue inequality $0<|F(x, y)| \leq h$,  improving an earlier result of Siegel and of Akhtari, Saradha \& Sharma.
\end{abstract}

\maketitle

\section{Introduction}

Let $F(x,y)$ be a binary form  with integer coefficients, having at least three pairwise non-proportional linear factors in its factorization over $\mathbb{C}$, and let $h$ be
a non-zero integer. Thue \cite{Thu1} proved that the equation 
$$F(x,y)=h$$
has only a finite number of integer solutions.  Such equations are called Thue equations. The problem of bounding the number of solutions of such equations, when the discriminant of $F$ has large absolute value, has garnered significant interest.
See, for example, \cite{Akh0}, \cite{Akhs}--\cite{Gy2}, 
\cite{SS1}--\cite{Ste} \& \cite{Wak1}. 

Consider a {\it diagonalizable} binary form $F(x,y)\in \mathbb Z[x,y]$ given by
$$F(x,y)=(\alpha x+\beta y)^r-(\gamma x+\delta y)^r$$
where the constants $\alpha,\beta,\gamma$ and $\delta$ satisfy
$$j=\alpha \delta-\beta\gamma \neq 0.$$
It turns out (see \cite[Lemma 4.1]{Vou}) that either 
$\alpha,\beta,\gamma,\delta\in\mathbb{Q}$ or $[\mathbb{Q}(\beta/\alpha):\mathbb{Q}]=2$, $\delta/\gamma$ is the algebraic conjugate of $\beta/\alpha$ over $\mathbb{Q}$, $\alpha^r\in \mathbb{Q}(\beta/\alpha)$ and $-\gamma^r$ is the algebraic conjugate of $\alpha^r$ over $\mathbb{Q}$.
Hence
$$uv=(\alpha x+\beta y)(\gamma x+\delta y)=\chi(Ax^2+Bxy+Cy^2)$$
for some $A,B,C\in \mathbb Z$ and a constant $\chi,$ where $\chi^r\in\mathbb{Q}$. Let
$$D=D(F)=B^2-4AC.$$
Then
$$j^2=\chi^2D.$$
Thus $D\neq 0.$ 
 We denote the discriminant of $F(x,y)$ by $\Delta=\Delta(F).$  
	 Put
\begin{equation}\label{Delta}
\Delta'=\frac{|\Delta|}{2^{r^2-r}r^rh^{2r-2}}.
\end{equation}

Let $h\in \mathbb{N}$ and consider the Thue inequality
\begin{equation}\label{thue}
0<|F(x,y)|\leq h.
\end{equation}
We are interested in counting the number of  non-zero tuples $(x,y)\in \mathbb Z \times \mathbb Z$ satisfying this inequality. Such a solution $(x,y)$ is said to be primitive if $\gcd(x,y)=1.$ We count $(x,y)$ and $(-x,-y)$ as one solution. Let $N_F(h)$ denote the number of primitive solutions to the inequality \eqref{thue}.

In 1970, Siegel \cite{Sie} proved the following theorem.
\begin{theorem}\label{Siegel}
Suppose that
$$\Delta'>(r^4h)^{c_\ell r^{2-\ell}}$$
where $r\geq 6-\ell, \ell=1,2,3,$
$$c_1=45+\frac{593}{913}, c_2=6+\frac{134}{4583}\ and \ c_3=75+\frac{156}{167}.$$
Then
$$ N_F(h)\leq 
\begin{cases}
2\ell r\ \textrm{if} \ D<0,\\
4\ell\ \ \textrm{if} \ D>0,\ r\ \textrm{is even and}\ F\ \textrm{is indefinite},\\
2\ell\ \ \textrm{if} \ D>0,\ r\ \textrm{is odd and}\ F\ \textrm{is indefinite},\\
1\  \textrm{if} \ D>0,\ \textrm{and}\ F\ \textrm{is definite}.\\
\end{cases}
$$
\end{theorem}
In 2018, Akhtari, Saradha \& Sharma \cite{AkSaSh} improved  the above result as follows.
See \cite[Theorems 1.3 and 1.4]{AkSaSh}.
\begin{theorem}\label{AkSS}
(i)
Let $r\geq 6$ and
\begin{equation}\label{i}
\Delta'\geq r^{13r^2(r-1)/(r^2-5r-2)}h^{4(r-1)(r^2-r+2)/(r^2-5r-2)}.
\end{equation}
Then
$$N_F(h)\leq 
\begin{cases}
2r+1\ \textrm{if} \ D<0,\\
5\ \textrm{if}\ D>0, r\ \textrm{is even and}\ F\ \textrm{is indefinite},\\
3\ \textrm{if}\ D>0, r\ \textrm{is odd and} \ F\ \textrm{is indefinite},\\
1\ \textrm{if}\ D<0\ \textrm{and}\ F\ \textrm{is definite}.\\
\end{cases}
$$
(ii)
Let $r\geq 5$ and
\begin{equation}\label{ii}
\Delta'\geq r^{7r^2(r-1)/((r-1)^{m-1}-2r-1)}h^{(r-1)(r^2+r+2)/((r-1)^{m-1}-2r-1)}.
\end{equation}
Then for any $m\geq 3$ we have
$$N_F(h)\leq 
\begin{cases}
rm\ \textrm{if} \ D<0,\\
2m\ \textrm{if}\ D>0, r\ \textrm{is even and}\ F\ \textrm{is indefinite},\\
m\ \textrm{if}\ D>0, r\ \textrm{is odd and}\ F\ \textrm{is indefinite},\\
1\ \textrm{if}\ D<0\ \textrm{and}\ F\ \textrm{is definite}.\\
\end{cases}
$$
\end{theorem}
We refer to \cite[Table 1]{AkSaSh} for a comparison of Theorem \ref{AkSS} with Theorem \ref{Siegel}.  In short, Theorem \ref{AkSS}(ii) is better than Theorem \ref{Siegel} as far as the lower bound for $\Delta'$ is concerned and give the same upper bound for $N_F(h)$
by taking $m=2\ell, \ell=2,3.$ For the case $\ell=1$ Theorem \ref{AkSS} is better than
Theorem \ref{Siegel} with respect to $\Delta'.$ But in the estimate for $N_F(h)$, we overshot by 1.  This discrepancy seems to be coming out of a {\it technical glitch} while counting the number of solutions of \eqref{inequality} that are related to an $r$-th root of unity. See Section 3 for details. In this paper, we correct this discrepancy and improve Theorem \ref{AkSS} as follows.
 \begin{theorem}\label{Main}
Let $r\geq 7$ and
\begin{equation}\label{ii}
\Delta'\geq r^{13r^2(r-1)/(r^2-5r-2)}h^{4(r-1)(r^2-r+2)/(r^2-5r-2)}.
\end{equation}
Then
$$N_F(h)\leq 
\begin{cases}
2r\ \textrm{if} \ D<0,\\
4\ \textrm{if}\ D>0, r\ \textrm{is even and}\ F\ \textrm{is indefinite},\\
2\ \textrm{if}\ D>0, r\ \textrm{is odd and} \ F\ \textrm{is indefinite},\\
1\ \textrm{if}\ D<0\ \textrm{and}\ F\ \textrm{is definite}.\\
\end{cases}
$$
\end{theorem}
Note that the bound for $\Delta'$ in \eqref{ii} is better than that of Siegel's in Theorem \ref{Siegel} when $\ell=1.$

\section{Notation and preliminaries}
Let $(x,y)\in \mathbb Z^2$ be a generic primitive solution of \eqref{thue} i.e., the inequality
\begin{equation*}
0<|F(x,y)|\leq h.
\end{equation*}
If $f$ is any function of $(x,y)$ then we write
$$f=f(x,y).$$
While enumerating the solutions of \eqref{thue} as $(x_0,y_0),(x_1,y_1),\ldots$ we denote by
$$f_i=f(x_i,y_i), i\geq 0.$$ 
We follow this notation throughout the paper without further mention.
Define the functions $u,v,\xi,\eta,\mu,Z,\zeta$ of $(x,y)$ as follows.
$$u=\alpha x+\beta y, v=\gamma x+\delta y, \xi=u^r,\ \eta=v^r.$$
From \cite[Lemmata 3.1 \& 3.2]{AkSaSh}, it follows that
when $D<0, u$ and $v$ are complex conjugates, giving $|u|=|v|.$ When $D>0, \xi$ and
$\eta$ are real with $u$ and $v$ as algebraic conjugates. 
Note that
$$F(x,y)=\xi-\eta=u^r-v^r.$$
Define
$$\mu=\frac{\eta}{\xi}=\frac{v^r}{u^r},\ Z=\max(|u|,|v|),\  \zeta=\frac{|F|}{Z^r}.$$
Note that when $D>0, \mu$ is real and
\begin{equation}\label{zeta}
\zeta=
\begin{cases}
1-\mu\ \textrm{if}\ 0<\mu<1,\\
 1-\mu^{-1}\ \textrm{if} \ \mu>1,\\
 1+|\mu|\ \textrm{if} \ -1<\mu<0,\\
 1+|\mu^{-1}|\ \textrm{if}\ \mu<-1.
\end{cases}
\end{equation}
\section{Solutions {\it related to} an $r-$th root of unity}
Write
$$F=\xi-\eta=\prod_{k=1}^r(u-v e^{\frac{2\pi i k}{r}}).$$
Let $\omega$ be an $r$-th root of unity. We say that the solution $(x,y)$ is {\it related to} 
$\omega$ if
$$|u-v\omega|=\min_{1\leq k\leq r}|u-v e^{\frac{2\pi i k}{r}}|.$$
Assume that
$$\frac{u}{v}=\bigg|\frac{u}{v}\bigg|e^{i\theta}$$
where
$$\frac{2(h-1)\pi}{r}\leq \theta<\frac{2h\pi}{r}$$
for some integer $h$ with $1\leq h< r$ i.e. $u/v$ lies in the arc subtended by the rays 
passing through the two $r$-th roots of unity viz.,  $e^{\frac{2(h-1)\pi}{r}}$ and $e^{\frac{2h\pi}{r}}.$ The perpendicular bisector of the chord joining these two roots of unity is $z=e^{\frac{h\pi}{r}}.$ So if $u/v$ lies {\it below} this line, then $\omega=
e^{\frac{2(h-1)\pi}{r}}.$ If $u/v$ lies {\it above} this line, then $\omega=
e^{\frac{2h\pi}{r}}.$ If $u/v$ lies on this line then it is at equal distance from both the roots of unity. In this case, as a convention we take $\omega=
e^{\frac{2(h-1)\pi}{r}}.$ Thus we see that every solution is {\it uniquely} related to an $r$-th root of unity.
We denote by $S$ the set of all primitive solutions of \eqref{thue} and by $S_\omega$ the set of all primitive solutions of \eqref{thue} that are related to $\omega.$ Then
$$S=\cup S_\omega$$
where $\omega $ ranges over all the $r$-th roots of unity. Thus it is enough to estimate $|S_\omega|.$ In the following lemmas we restrict to the set $S_\omega.$ Let us enumerate all the solutions in $S_\omega$ as $(x_1,y_1),(x_2,y_2),\ldots$ so that $\cdots \leq \zeta_3\leq \zeta_2\leq \zeta_1.$ 
Here we re-do many of the lemmas from \cite[Section 5]{AkSaSh}. We point out that many of these lemmas are proved under the assumption that the relevant $\zeta<1.$ Below we 
avoid this assumption. This is the {\it technical glitch} referred to in the Introduction.
\begin{lemma}\label{realmu}
Let $D>0.$
Suppose $(x,y)\in S_\omega.$ Then
$$\frac{u}{v}=
|\mu^{-1}|^{1/r}e^{\pi i\epsilon/r}\omega
$$
where
$$\epsilon=
\begin{cases}
0\ if\ \mu>0\\
1\ if \ \mu<0.
\end{cases} 
$$
\end{lemma}

\begin{proof}
As noted earlier, when $D>0,\mu^{-1}$ is real. So we write
$$\frac {u^r}{v^r}=|\mu^{-1}|e^{\pi i\epsilon }.$$
Then
\begin{equation}\label{mu0}
\frac{u}{v}=|\mu^{-1}|^{1/r} e^{\pi i(\epsilon+2h)/r}
\end{equation}
for some integer $h$ with $0\leq h<r.$ Since $(x,y)\in S_\omega,$ we have
\begin{equation}\label{mu}
\bigg||\mu^{-1}|^{1/r}e^{\pi i\epsilon/r} -\omega e^{-2\pi i h/r}\bigg|=
\min_{0\leq k<r}\bigg||\mu^{-1}|^{1/r}e^{\pi i\epsilon/r} -e^{2\pi i(k-h)/r}\bigg|. 
\end{equation}
As $k$ varies from $0$ to $r-1, k-h$ varies over a complete residue system $\pmod r.$ 
So
$$\min_{0\leq k<r}\bigg||\mu^{-1}|^{1/r}e^{\pi i \epsilon/r} -e^{2\pi i(k-h)/r}\bigg|
=\min_{0\leq k<r}\bigg||\mu^{-1}|^{1/r}e^{\pi i\epsilon/r} -e^{2\pi  i k/r}\bigg|. $$ 
Now $|\mu^{-1}|^{1/r}e^{i\pi \epsilon/r}$ is a complex number in the first quadrant with argument $\pi \epsilon/r.$ So by our convention, the $r$-th root of unity nearest to it is 1. Hence we conclude in \eqref{mu} that
$\omega e^{- 2\pi ih/r}=1$ giving
$$ e^{2\pi ih/r}=\omega$$
which proves the result.

\end{proof}

\begin{lemma}\label{AllD}
Let $(x,y)\in S_\omega.$ Also assume that $\epsilon=0$ if $D>0.$ Then
$$\bigg|\frac{u}{v}-\omega\bigg|\leq \frac{Z}{|v|}\zeta .
$$
\end{lemma}
\begin{proof}
Let $D<0.$ Then from \cite[Lemma 5.5 ,(37)]{AkSaSh}, we have for $r\geq 3,$
$$\bigg|\frac{u}{v}-\omega\bigg|\leq \frac{\pi}{2r}\zeta\leq \zeta\leq \frac{Z}{|v|}\zeta.
$$
Now we take $D>0.$ Since $\epsilon=0,$ by Lemma \eqref{realmu}, we get
$$\frac{u}{v}=|\mu^{-1}|^{1/r}\omega.$$
Observe that
$$\begin{cases}
0<|\mu^{-1}|^{1/r}<1\ \textrm{ if } \ |u|\leq |v|\\
|\mu^{-1}|^{1/r}>1\ \textrm{ if } \ |u|> |v|.
\end{cases}
$$ 
So it follows from Lemma \ref{realmu}, that
$$
\bigg|\frac{u}{v}-\omega\bigg|=\bigg||\mu^{-1}|^{1/r}-1\bigg|=
\begin{cases}
1-|\mu^{-1}|^{1/r}\ \textrm{ if } \ |u|\leq |v|\\
|\mu^{-1}|^{1/r}-1\ \textrm{ if } \ |u|> |v|.
\end{cases}
$$
Let $|u|\leq |v|.$ Then
$$
\bigg|\frac{u}{v}-\omega\bigg|\leq 1-|\mu^{-1}|=\frac{|F|}{Z^r}=\zeta.$$
Let $|u|>|v|.$ Then
\begin{align*}
\bigg|\frac{u}{v}-\omega\bigg|&=
|\mu^{-1}|^{1/r}-1=|\mu^{-1}|^{1/r}(1-|\mu|^{1/r})\\
&\leq |\mu^{-1}|^{1/r}(1-|\mu|)= \frac{|u|}{|v|}\bigg(1-\frac{|v|^r}{|u|^r}\bigg)\\
&=\frac{|u|}{|v|}\frac{|F|}{Z^r}=\frac{Z}{|v|}\zeta.
\end{align*}


The proof of the lemma is complete.

\end{proof}
\begin{lemma}\label{ZZ*}
Let $(x,y)\neq (x_*,y_*)$ be two primitive solutions of \eqref{thue} in $S_\omega$ with $\zeta_*\leq \zeta.$  Then
$$Z_*\geq \frac{|j|}{2h^{1/r}}.$$
\end{lemma}
\begin{proof}
Consider
\begin{equation}\label{xy*}
uv_*-u_*v=(\alpha \delta-\beta\gamma)(xy_*-yx_*)=j(xy_*-yx_*)\neq 0.
\end{equation}
Hence
\begin{equation}\label{Z*}
|j|\leq |uv_*|+|u_*v|\leq 2Z Z_*.
\end{equation}
Thus
$$Z_*\geq \frac{|j|\zeta^{1/r}}{2|F|^{1/r}}$$
which gives the assertion if $\zeta\geq 1,$ since $|F|\leq h.$ So let us assume that $\zeta<1.$ 
If $D>0,$ then we have
$\mu >0, \zeta_*<1$ and $\mu_*>0$ so that Lemma \ref{AllD} is applicable.
Thus from \eqref{xy*} and Lemma \ref{AllD}, we get
\begin{align*}
|j|  &\leq |vv_*|\bigg|\frac{u}{v}-\frac{u_*}{v_*}\bigg|\leq |vv_*|\bigg(\bigg|\frac{u}{v}-\omega\bigg|+\bigg|\frac{u_*}{v_*}-\omega\bigg|\bigg)\\
& \leq |vv_*|\bigg(\frac{Z\zeta}{|v|}+\frac{Z_*\zeta_*}{|v_*|}\bigg)\leq 2ZZ_*\zeta.
\end{align*}
Hence 
\begin{equation}\label{1/r}
|j|\leq 2h^{1/r}\zeta^{(r-1)/r}Z_*\leq 2h^{1/r}Z_*
\end{equation}
which proves the assertion for $\zeta<1.$ 
\end{proof}
\noindent
{\bf Note.}

From \eqref{Z*},
$$|j|\leq 2h^{2/r}\zeta^{-1/r}\zeta_*^{-1/r}\leq 2h^{2/r}\zeta_*^{-2/r}$$
giving
$$\zeta_*\leq 2^{r/2}h |j|^{-r/2}.$$
Thus if $|j|>2^{1+2\nu/r} h^{2/r},$ then $\zeta_*<2^{-\nu}.$
In particular, if $|j|>2h^{2/r},$ then $\zeta_*<1.$ 
As before, let us arrange the solutions in $S_\omega$ as
$(x_1,y_1),(x_2,y_2),\ldots,$ so that
$$\zeta_1\geq \zeta_2\geq \zeta_3\geq\cdots.$$
By our observation above, we see that 
\begin{equation}\label{less1}
\zeta_i<1\ \textrm{for} \ i \geq 2\ \textrm{if} \ |j|>2h^{2/r}.
\end{equation}
We will use this arrangement of the solutions in $S_\omega$ from now on.
\begin{lemma}\label{iteration}
Let $(x_1,y_1),(x_2,y_2),\ldots, (x_t,y_t)$ be in $S_\omega$ with $t\geq 3.$
Assume that
$$|j|>2^{1+(r-2)/(r(R(t-1)-1))}h^{2/r}.$$
Then
$$\zeta_{t-1}<1/2.$$
\end{lemma}
\begin{proof}
Applying Lemma \ref{ZZ*} with $(x,y)=(x_1,y_1)$ and $(x_*,y_*)=(x_2,y_2)$ we get
$$\zeta_2\leq \frac{2^rh^2}{|j|^r}.$$
Thus
$$\zeta_2\leq H^r$$
where
$$H=2h^{2/r}|j|^{-1}.$$
Again applying Lemma \ref{ZZ*} with $(x,y)=(x_2,y_2)$ and $(x_*,y_*)=(x_3,y_3)$ we get
$$\zeta_3\leq H^{r}\zeta_2^{r-1}\leq H^{r(1+(r-1))}.$$
Proceeding inductively, we obtain that
$$\zeta_{t-1}\leq H^{r(1+(r-1)+\cdots+ (r-1)^{t-3})}.$$
Thus $\zeta_{t-1}<1/2$ if
$$|j|^{r(R(t-1)-1)/(r-2)}>2^{1+r(R(t-1)-1)/(r-2)}h^{2(R(t-1)-1)/(r-2)}.$$
Hence $\zeta_{t-1}<1/2$ if
$$|j|>2^{1+(r-2)/(r(R(t-1)-1))}h^{2/r} .$$

\end{proof}

\begin{lemma}\label{gap-principle}
Let $(x_{i-1},y_{i-1}),(x_i,y_i)$ be in $S_\omega$ with $i\geq 2.$ 
Assume that $|j|>2h^{2/r}$ if $D>0.$
Then
$$Z_i\geq \frac{|j|}{2h}Z_{i-1}^{r-1}.$$
\end{lemma}
\begin{proof}
We always have 
$$\zeta_i\leq \zeta_{i-1}.$$
Note that
$$|j| \leq |u_{i-1}v_i-u_iv_{i-1}|.$$
So it is enough to show that
\begin{equation}\label{ZZ}
|u_{i-1}v_i-u_iv_{i-1}|\leq 2Z_{i-1}Z_i\zeta_{i-1}
\end{equation}
since $\zeta_{i-1}\leq h/Z_{i-1}^r.$
First let $D<0.$ Then by Lemma \ref{AllD}
\begin{align*}
|u_{i-1}v_i-u_iv_{i-1}| & \leq |v_{i-1}v_i|\bigg(\bigg|\frac{u_{i-1}}{v_{i-1}}-\omega\bigg|+\bigg|\frac{u_{i}}{v_{i}}-\omega\bigg|\bigg)\\
& \leq Z_{i-1}Z_i(\zeta_{i-1}+\zeta_i)\leq 2Z_{i-1}Z_i\zeta_{i-1}.
\end{align*}

Next let $D>0.$ Then by \eqref{less1}, $\zeta_i<1.$ By Lemma \ref{realmu} we have
$$\frac{u_{i-1}}{v_{i-1}}=|\mu_{i-1}^{-1}|^{1/r}e^{\pi i\epsilon}\omega;\ \
\frac{u_i}{v_i}=|\mu_{i}^{-1}|^{1/r}\omega.$$
So
\begin{equation}\label{four}
 |u_{i-1}v_i-u_iv_{i-1}|=
\begin{cases}
(i) \ |u_{i-1}u_i|\ \bigg||\mu_{i}|^{1/r}-|\mu_{i-1}|^{1/r}e^{-\pi i\epsilon/r}\bigg|\\
(ii) \ |u_{i-1}v_i|\ \bigg|1-|\mu_{i}^{-1}|^{1/r}|\mu_{i-1}|^{1/r}e^{-\pi i\epsilon/r}\bigg|\\
(iii) \ |v_{i-1}u_i|\ \bigg|1-|\mu_{i}|^{1/r}|\mu_{i-1}^{-1}|^{1/r}e^{\pi i\epsilon/r}\bigg|\\
(iv) \ |v_{i-1}v_i|\ \bigg||\mu_{i}^{-1}|^{1/r}-|\mu_{i-1}^{-1}|^{1/r}e^{\pi i\epsilon/r}\bigg|.
\end{cases}
\end{equation}

First we deal with the case $\epsilon=0.$ Then $\mu_{i-1}>0.$
Since $\zeta_i<1$, by \eqref{zeta} we get $\mu_i>0.$ We need to consider 8 cases depending on the signs and the values of $\mu_{i-1}$ and $\mu_{i}.$ In each case we show that \eqref{ZZ}
holds.
\vskip 1mm
\noindent
{\it Case 1.} Let $0<\mu_i,\mu_{i-1}<1.$\\
Then
$$0<1-\mu_{i-1}^{1/r}<1-\mu_{i-1}=\zeta_{i-1};\ 0<1-\mu_{i}^{1/r}<1-\mu_i=\zeta_i.$$
By \eqref{four}(i),
\begin{align*}
|u_{i-1}v_i-u_iv_{i-1}|& =|u_{i-1}u_i|\ \bigg|(1-|\mu_{i}|^{1/r})-(1-|\mu_{i-1}|^{1/r})\bigg|\\
&\leq |u_{i-1}u_i|(\zeta_{i-1}+\zeta_i)\\
&\leq 2Z_{i-1}Z_i\zeta_{i-1}.
\end{align*}
\vskip 1mm
\noindent
{\it Case 2.} Let $0<\mu_i<1, \mu_{i-1}>1.$\\
Then
$$0<1-\mu_{i-1}^{-1/r}<1-\mu_{i-1}^{-1}=\zeta_{i-1};\ 0<1-\mu_{i}^{1/r}<1-\mu_i=\zeta_i$$
and
$$0<1-\mu_{i-1}^{-1/r}\mu_i^{1/r}<1-\mu_{i-1}^{-1}\mu_i=1-(1-\zeta_{i-1})(1-\zeta_i)<2\zeta_{i-1}.$$
By \eqref{four}(iii),
\begin{align*}
|u_{i-1}v_i-u_iv_{i-1}|& =|v_{i-1}u_i|\ \bigg|1-\mu_{i}^{1/r}|\mu_{i-1}|^{1/r}\bigg|\\
&\leq 2|v_{i-1}u_i|\zeta_{i-1}\\
&\leq 2Z_{i-1}Z_i\zeta_{i-1}.
\end{align*}



\vskip 1mm
\noindent
{\it Case 3.} Let $\mu_i>1,0<\mu_{i-1}<1.$\\
Then
$$\zeta_{i-1}=1-|\mu_{i-1}|;\ 0<1-\mu_{i}^{-1/r}<1-\mu_i^{-1}=\zeta_i$$

By \eqref{four}(ii),
\begin{align*}
|u_{i-1}v_i-u_iv_{i-1}|& =|u_{i-1}v_i|\ \bigg|1-|\mu_{i}|^{-1/r}|\mu_{i-1}|^{1/r}\bigg|\\
&\leq |u_{i-1}v_i|(1-|\mu_{i-1}\mu_i^{-1}|)\\
&\leq |u_{i-1}v_i|(1-|\mu_{i-1}|(1-\zeta_i))\\
&\leq |u_{i-1}v_i|(1+\zeta_i)\leq 2Z_{i-1}Z_i\zeta_{i-1}.
\end{align*}

\vskip 1mm
\noindent
{\it Case 4
.} Let $\mu_i>1,\mu_{i-1}>1.$\\
Then
$$0<1-\mu_{i-1}^{-1/r}<1-\mu_{i-1}^{-1}=\zeta_{i-1};\ 0<1-\mu_{i}^{-1/r}<1-\mu_i^{-1}=\zeta_i$$

By \eqref{four}(iv),
\begin{align*}
|u_{i-1}v_i-u_iv_{i-1}|& =|v_{i-1}v_i|\ \bigg||\mu_{i}|^{-1/r}-|\mu_{i-1}|^{1/r}\bigg|\\
&\leq |v_{i-1}v_i|(|1-\mu_{i-1}^{-1}|+|1-\mu_i^{-1}|)\\
&\leq |v_{i-1}v_i|(\zeta_{i-1}+\zeta_i)\\
&\leq 2|v_{i-1}v_i|\zeta_{i-1}\leq 2Z_{i-1}Z_i\zeta_{i-1}.
\end{align*}
This completes all the cases when $\epsilon=0.$
\vskip 2mm
Next we take $\epsilon=1.$ Then $\mu_{i-1}<0$ and so $\zeta_{i-1}>1.$ Recall that $\mu_i>0.$ There are four cases to consider, viz.,
$$(a) \ 0<\mu_i\leq 1, -1\leq \mu_{i-1}<0 \ (b)\ 0<\mu_i\leq 1, \mu_{i-1}<-1;$$
$$(c) \ \mu_i> 1, -1\leq \mu_{i-1}<0 \ 
(d) \mu_i> 1, \mu_{i-1}<-1.$$ 
Using  (i), (iii), (ii) and (iv) respectively for (a),(b),(c) and (d), we find that
$$|u_{i-1}v_i-u_iv_{i-1}|\leq 2Z_{i-1}Z_i<2 Z_{i-1}Z_i\zeta_{i-1}$$
since $\zeta_{i-1}>1.$

\end{proof}
\section{Condition on $|j|$}
Using Pad\'e approximation, in \cite[Lemma 7.3]{AkSaSh} certain algebraic numbers were constructed. This construction was used along with symbolic computation and Lemma \ref{gap-principle} to get an iterative bound in \cite[Lemma 8.1]{AkSaSh} which we state below. As the Lemmas in Section 3 are valid without any condition on the $\zeta'$s, we are able to state this lemma for $S_\omega$ rather than $S_\omega'$ (see \cite[(40)]{AkSaSh}). We shall also restrict to $|S_\omega|=3.$   
\begin{lemma}\label{comp}
Let $r\geq 7.$ Assume that $|S_\omega|= 3.$ Suppose that
\begin{equation}\label{j}
|j|\geq 2r^{i_7/r}h^{i_8/r}
\end{equation}
with
\begin{equation}\label{i7}
i_7=\frac{13r^2}{r^2-5r-2}, \ i_8=\frac{2(3r-1)(r-2)}{r^2-5r-2}.
\end{equation}
Then for every integer $n \geq 1,$ we have
\begin{equation}\label{Zk}
Z_3\geq \frac{Z_{2}^{nr}}{2^{n+4}r^{(3nr+2)/(r-2)}|j|^{(nr+2)/(r-2)}h^{2n+1}}.
\end{equation}
\end{lemma}
Note that the exponent of $Z_2$ in the above estimate is $nr$ while in Lemma 8.1 of \cite{AkSaSh}, it was $(n+1)r-1.$ This weaker estimate is sufficient for our purpose.
\begin{proof}
We follow the proof of Lemma 8.1 of \cite{AkSaSh}. We give details wherever necessary. 
Let us denote the solutions in $S_{\omega}$ as
$$(x_1,y_1),(x_2,y_2),(x_3,y_3)\ \textrm{ with }\ \zeta_1\geq \zeta_2\geq \zeta_3.$$ 
Then
\begin{equation}\label{ZZ}
Z_3\geq \frac{|j|}{2h} Z_2^{r-1},  Z_2\geq \frac{|j|}{2h^{1/r}}
\end{equation}
giving
$$Z_3\geq \frac{|j|^r}{2^r h^{2-1/r}}.$$
Now we follow the argument as in the proof of Theorem 1.3 of \cite{AkSaSh} to get
\begin{equation}\label{761}
Z_3\geq \frac{Z_2^{r(n+1-g)-1+g}}{2^{n+4}r^{(r(2g+3n)+2)/(r-2)}|j|^{(r(g+n)+2)/(r-2)} h^{2n+1-g}},n\geq 1, g\in \{0,1\}.
\end{equation}
if
\begin{equation}\label{762}
Z_3^{r-1}\geq 2^{3n+4}r^{(r(2g+3n)+2)/(r-2)}|j|^{(r(g+n)+2)/(r-2)} h Z_2^{nr+1-g}.
\end{equation}
For a given integer $n\geq 1$ and $g \in \{0,1\}$ let $a_i=a_i(n,r,g), 1\leq i\leq 5,$ we say that property $P[a_1,a_2,a_3,a_4,a_5]$ holds if 
\begin{equation}\label{Z1Z2}
Z_3\geq \frac{Z_2^{a_1}}{2^{a_2}r^{a_3}|j|^{a_4}h^{a_5}}.
\end{equation} 
Let us now assume that $P[a_1,a_2,a_3,a_4,a_5]$ holds with $a_2+a_4\geq 0.$
Then \eqref{762} is valid if
\begin{equation*}
Z_2^{a_1(r-1)-nr-1+g}\geq 
2^{a_2(r-1)+3n+4}r^{a_3(r-1)+(r(2g+3n)+2)/(r-2)}|j|^{a_4(r-1)+(r(g+n)+2)/(r-2)}
h^{a_5(r-1)+1}.
\end{equation*}
By \eqref{ZZ}, the above inequality is valid if
\begin{equation}\label{j1}
|j|^{A_1}\geq 2^{A_1+a_2(r-1)+3n+4}r^{a_3(r-1)+(r(2g+3n)+2)/(r-2)}|j|^{a_4(r-1)+(r(g+n)+2)/(r-2)}
h^{A_1/r+a_5(r-1)+1}
\end{equation}
where 
$$A_1=a_1(r-1)-nr-1+g.$$
We set
\begin{align*}
B_1=&A_1-a_4(r-1)-(r(g+n)+2)/(r-2)\\
 B_2=&A_1+a_2(r-1)+3n+4\\
B_3=&a_3(r-1)+(r(2g+3n)+2)/(r-2)\\
B_4=& A_1/r+a_5(r-1)+1.
\end{align*}
We implement the induction procedure given in \cite[Section 8]{AkSaSh} with the above values of $B_1,\cdots, B_4$ and $a_1,\cdots, a_5$ as already given in the procedure.
If the conditions 
$$(i) \ A_1>0\quad (ii)\ B_1>0\quad (iii)\ B_1\times \frac{13r^2}{r^2-5r-2}\geq r(B_3+(B_2-B_1)/2)$$
$$(iv)\ B_1\times \frac{2(3r-1)(r-2)}{r^2-5r-2}\geq r B_4$$
are satisfied then \eqref{j1} holds and hence \eqref{761} holds.
To begin with, by \eqref{Z1Z2}, $P[r-1,1,0,-1,1]$ holds.
Let $(n,g)=(1,0).$ The conditions (i)--(iv) are satisfied for $a_1=r-1,a_2=1,a_3=0,a_4=-1,a_5=1$.
If $\Sigma_{1,0}\neq 0$  (see \cite[Section 7]{AkSaSh} for the definition of $\Sigma_{n,g}$), we get that
$P[ 2r-1,5,(3r+2)/(r-2),(r+2)/(r-2),3]$ is valid. Thus \eqref{Zk} is valid with $n=1$ since
$$Z_2^{r-1}\geq 1$$
by \eqref{j},\eqref{i7} and \eqref{ZZ}. 
 If $\Sigma_{1,0}= 0$, then $\Sigma_{1,1},\Sigma_{2,1}$ are non-zero. 
Fixing $(n,g)=(1,1)$
and using the same parameters as in \cite[Lemma 8.1]{AkSaSh}, we find that
$P[ r,5,(5r+2)/(r-2),(2r+2)/(r-2),2]$ is valid.  Using this and taking $(n,g)=(2,1)$, we get that
$$Z_3\geq \frac{Z_2^{2r}}{2^6r^{(8r+2)/(r-2)}|j|^{(3r+2)/(r-2)}h^4}.$$
Thus \eqref{Zk} is valid with $n=1$ provided 
\begin{equation}\label{673}
\frac{Z_2^{r}}{2r^{5r/(r-2)}|j|^{2r/(r-2)}h}\geq 1.
\end{equation}
This is again satisfied by \eqref{j},\eqref{i7} and \eqref{ZZ}.
Now we proceed by induction as in \cite[Lemma 8.1]{AkSaSh} to complete the lemma.
\end{proof}

\section{Proof of Theorem \ref{Main}}
Assume that \eqref{i} holds.  Since
	\[
 	 \Delta = (-1)^{\frac{(r-1)(r+2)}{2}} r^r j^{r(r-1)},
 	 \]
we get \eqref{j}.
Now,
by \cite[Theorem 1.3]{AkSaSh}, we may assume that
$$N_F(h)=|S|=2r+1.$$
Our arguments are similar to Theorem 1.3 of \cite{AkSaSh}. We give the details wherever it is necessary.
By our assumption on $N_F(h),$ there exists some $\omega,$ say $\omega_1$ with
$$|S_{\omega_1}|\geq 3.$$
Suppose there is another $\omega_2$ with $|S_{\omega_2}|\geq 3.$ Then the solution
$(x_0,y_0)$ with $\zeta_0$ largest can belong to at most only one of these two sets. Hence there exists a set $S_\omega, \omega=\omega_1$ or $\omega_2$ such that
$S_\omega$ does not contain $(x_0,y_0)$ and $|S_\omega|\geq 3.$ 
Then $S_\omega=S_{\omega"}$ where $S_{\omega"}$ is as in the proof of Theorem 1.3 in \cite{AkSaSh}.
So we may follow the argument therein to get a contradiction.
Hence there exists exactly one $\omega$, say $\omega_1$ such that
$$|S_{\omega_1}|=3\ \textrm{and} \ |S_\omega|=2\ \textrm{for} \ \omega\neq \omega_1.$$  
Now by Lemma \ref{comp} applied to $S_{\omega_1}$ we get, with the same notation as in there, that
\begin{equation}\label{Zk1}
Z_3\geq \frac{Z_{2}^{nr}}{2^{n+4}r^{(3nr+2)/(r-2)}|j|^{(nr+2)/(r-2)}h^{2n+1}}
\end{equation}
holds for all $n\geq 1$.
But the right hand side approaches infinity as $n$ approaches infinity, which gives the final contradiction.
\qed

\section*{Acknowledgement}
Saradha likes to thank the Indian National Science Academy for awarding the Senior Scientist fellowship under which this work was done. She also thanks DAE-Center for Excellence in Basic Sciences, Mumbai University for providing facilities to carry out this work. Divyum
acknowledges the support of the DST-SERB SRG Grant SRG/2021/000773 and the OPERA award of BITS Pilani.

\end{document}